\input btxmac.tex

\long\def\comment#1{\relax}


\magnification=\magstep1 
\hsize=15truecm 
\vsize=22.3truecm
\hoffset=.65truecm 
\voffset=.3truecm


\hyphenation{Null-stel-len-satz Pro-po-si-tion}


\font\tenbb=msbm10     
\font\sevenbb=msbm7
\newfam\bbfam \def\bb{\fam\bbfam\tenbb}  
\textfont\bbfam=\tenbb
\scriptfont\bbfam=\sevenbb

\font\tenam=msam10     
\font\sevenam=msam7
\newfam\amfam   
\textfont\amfam=\tenam
\scriptfont\amfam=\sevenam

\font\bn=cmb10

\font\titlefont=cmss12
\font\authorfont=cmr12
\let\titlerm=\elevenrm

\let\footnotefont=\eightandahalfrm
\let\dinky=\sevenrm


\def\glq{\hbox{,\kern-.7pt,\kern.3pt}}  
\chardef\colon="3A

\let\em=\it

\font\sevensl=cmsl8 at 7pt
\scriptfont\slfam=\sevensl

\let\titlefont=\twelverm  
\let\authorfont=\tenrm  
\let\sectheadfont=\titleit  

\long\def\suppress#1{\relax}


\def\secthead#1#2{\goodbreak\removelastskip\bigskip%
{\titlerm #1\enspace}{\sectheadfont #2}\par%
\nobreak\medskip\noindent}


\thickmuskip=4.5mu plus 4mu minus 3mu  
\mathsurround=0.5pt

\pretolerance=200  \hyphenpenalty=1000 \exhyphenpenalty=5000

\baselineskip=14pt \parskip=1.5pt \parindent=1.5em
\smallskipamount=4pt plus2pt minus1pt 
\medskipamount=7pt plus3pt minus2pt 
\bigskipamount=15pt plus5pt minus4pt 


\def\modulo#1#2{{\raise.15em\hbox{$#1$}/\!\lower.18em\hbox{$#2$}}}

\font\bigsy=cmsy10 at 12pt
\def\bigwr{\raise.255555pt\hbox{$ \wr \textfont2=\bigsy $}}

\font\smallex=cmex10 at 7pt
\def\smallsum{\lower1pt\hbox{$ \sum \textfont3=\smallex $}}

\def\pileftarrow#1{\buildrel{\,\pi_{#1}}\over\longleftarrow}
\def\namedleftarrow#1{\buildrel{\,#1}\over\longleftarrow}

\def\projlim{\lower6pt\hbox{$\buildrel\hbox{\tenrm lim}\over\leftarrow$}\ }

\let\originalsqrt\sqrt
\def\varsqrt #1#2{{\setbox0=\hbox{$#1\originalsqrt{#2\hskip2pt}$}%
     \dimen0=\ht0 \advance\dimen0-.5pt
     \dimen2=\dimen0
     \advance\dimen2 -1.5pt
     \ifdim \dimen2 >10pt \advance\dimen2-.5pt \fi
     \ifdim \dimen2 >15pt \advance\dimen2 -1pt \fi
     \box0 \kern-.4pt
     \vrule width.4pt height\dimen0 depth -\dimen2 \kern.8pt
  }}


\chardef\at="40  

\def\N{{\bb N}}     
\def\Z{{\bb Z}}     
\def\Zp{{\Z_{p}}} 
\def\Zx{{\Z[x]}} 
\def\Zpn{{\Z_{p^n}}} 
\def\Zpm{{\Z_{p^m}}} 
\def\Zptwo{{\Z_{p^2}}} 
\let\id\iota

\def\[{{\rm[}} \def\]{\/{\rm]}}  
\def\({{\rm(}} \def\){\/{\rm)}}  
\let\congr=\equiv    
\let\takeaway=\setminus  
\mathchardef\ideal="1945  
\mathchardef\eop="1903  

\newcount\codenumber \codenumber=`:
\loop\relax\if\codenumber>"FFF \advance\codenumber\by -"1000\repeat
\advance\codenumber by "3000
\mathcode`:=\codenumber


\long\def\pronounce#1#2\endpronounce {\ifdim\lastskip<\medskipamount%
  \removelastskip\penalty-50\medskip\fi
  \noindent{\bn#1.\enspace}{\sl#2\par}%
  \ifdim\lastskip<\medskipamount\removelastskip\penalty55\medskip\fi}

\long\def\rem#1#2\endrem {\ifdim\lastskip<\medskipamount%
  \removelastskip\medskip\penalty-350\fi
  \noindent{\bn#1.\ }{\rm#2\par}%
  \ifdim\lastskip<\smallskipamount\removelastskip\penalty-100\smallskip\fi}

\def\proof{\noindent {\it Proof.\enspace}}
\def\endproof{\ $\eop$\par%
   \ifdim\lastskip<\medskipamount
   \removelastskip\penalty-50\medskip\fi}


\def\firstpagehead{\tenrm To appear in J.~Number~Th. \hfill }
\def\runninghead{\hfill 
{\dinky  SYLOW GROUPS OF GROUPS OF POLYNOMIAL PERMUTATIONS} \hfill }
\headline={\ifnum\pageno>1\runninghead \else\firstpagehead \fi}

\output={
          \shipout\vbox{\makeheadline
                        \vskip.1truein
                        \pagebody
                        \makefootline }
          \advancepageno
          \ifnum\outputpenalty>-20000 \else\dosupereject\fi
 }


\line{\ \hfill\ }
\bigskip
\centerline{\titlefont Sylow $p$-groups of polynomial
permutations on the integers mod $p^n$}
\bigskip
\comment{
\centerline{\titlefont Sylow $p$-groups of polynomial permutations}
\smallskip
\centerline{\titlefont on the integers mod $p^n$}
\bigskip\medskip
}
\centerline{\authorfont Sophie Frisch\ \ and\ \ Daniel Krenn}
\footnote{}{\footnotefont Sophie Frisch is supported by the Austrian Science 
Fund (FWF), grant P23245-N18.}
\footnote{}{\footnotefont Daniel Krenn is supported by the Austrian Science 
Fund (FWF), project W1230 doctoral program ``Discrete Mathematics''. }
\bigskip\medskip

{\bf Abstract.}
We enumerate and describe the Sylow $p$-groups of the groups of 
polynomial permutations of the integers mod $p^n$ for $n\ge 1$ and 
of the pro-finite group which is ths projective limit of these groups.
MSC 2000: primary 20D20, secondary 11T06, 13M10, 11C08, 13F20, 20E18.
\medskip

\newif\ifkrenn \krennfalse
\newif\ifdissertationkrenn\dissertationkrennfalse

\ifkrenn\section{Introduction}\else
\secthead{1.}{Introduction}\fi
Fix a prime $p$ and let $n\in \N$. Every polynomial $f\in\Zx$
defines a function from $\Zpn= \modulo{\Z}{p^n\Z}$ to itself. If 
this function happens to be bijective, it is called a {\em polynomial
permutation} of $\Zpn$. The polynomial permutations of $\Zpn$
form a group $(G_n,\circ)$ with respect to composition. The 
order of this group has been known since at least 1921 
(Kempner \cite{Kem21PRS}) to be 
$$|G_2|=p!(p-1)^p p^p\quad\hbox{\rm\ and\ }
|G_n|=p!(p-1)^p p^p p^{\smallsum_{k=3}^n \beta(k)}\;\;
\hbox{\rm\ for\ } n\ge 3,$$
where $\beta(k)$ is the least $n$ such that $p^k$ divides $n!$,
but the structure of $(G_n,\circ)$ is elusive. 
(See, however, N\"obauer \cite{NoePFPR82} for some partial results). 
Since the order of $G_n$ is divisible by a high power of $(p-1)$ for 
large $p$, even the number of Sylow $p$-groups is not obvious. 

We will show that there are $(p-1)!(p-1)^{p-2}$ Sylow $p$-groups of $G_n$ and
describe these Sylow $p$-groups, see Theorem 5.1 and Corollary 5.2.

Some notation: $p$ is a fixed prime throughout.
A function $g\colon \Zpn\rightarrow \Zpn$ arising from a 
polynomial in $\Zpn[x]$ or, equivalently, from a polynomial in
$\Zx$, is called a {\em polynomial function} on $\Zpn$.
We denote by $(F_n,\circ)$ the monoid with respect to composition 
of polynomial functions on $\Zpn$.
By monoid, we mean semigroup with an identity element.
Let $(G_n,\circ)$ be the group of units of $(F_n,\circ)$, which
is the group of polynomial permutations of $\Zpn$.

Since every function induced by a polynomial preserves congruences
modulo ideals, there is a natural epimorphism mapping
polynomial functions on $\Z_{p^{n+1}}$ onto
polynomial functions on $\Zpn$, and we write it as 
$\pi_n\colon F_{n+1}\rightarrow F_n$. If $f$ is a polynomial in
$\Z[x]$ (or in $\Zpm[x]$ for $m\ge n$) we denote the polynomial
function on $\Zpn[x]$ induced by $f$ by $[f]_{p^n}$. 

The order of $F_n$ and that of $G_n$ have been determined by Kempner
\cite{Kem21PRS} in a rather complicated manner. His results
were cast into a simpler form by N\"obauer \cite{Noe55GRP} and
Keller and Olson \cite{KeOl68PF} among others.
Since then there have been many generalizations of the order
formulas to more general finite rings 
\cite{Ros75PRFR,Nec80PTPI,BrMu92Gal,Fri97pfr,Bhar97Pord,JGSZ06,J10npf}.
Also, polynomial
permutations in several variables (permutations of $(\Zpn)^k$ defined
by $k$-tuples of polynomials in $k$ variables) have been looked into
\cite{Fri95wpps,Chen96pfZni,Zh04pfpp,WZ07sos,WZ09pp2v,LJ09pfnv}.

\ifkrenn\section{Polynomial functions and permutations}\else
\secthead{2.}{Polynomial functions and permutations}%
\nobreak\fi
To put things in context, we recall some well-known facts, to be
found, among other places, in \cite{Kem21PRS,Noe55GRP,Car64fp,KeOl68PF}.
The reader familiar with polynomial functions on finite rings is 
encouraged to skip to section~3. \ifdissertationkrenn This section does not
contain new material but reviews the state of the art. \else  
Note that we do not claim anything in section~2 as new. \fi

\ifkrenn\begin{definition}\else\rem{Definition}\fi
For $p$ prime and $n\in\N$, let
$$\alpha_p(n)=\sum_{k=1}^{\infty}\left[{n\over{p^k}}\right]
\qquad\hbox{\rm and}\qquad
\beta_p(n)=\min\{m\mid \alpha_p(m)\ge n\}.$$
If $p$ is fixed, we just write $\alpha(n)$ and $\beta(n)$.
\ifkrenn\end{definition}\else\endrem\fi

\ifkrenn\begin{notation}\else\rem{Notation}\fi
For $k\in\N$, let $(x)_k=x(x-1)\ldots(x-k+1)$ and $(x)_0=1$. 
We denote $p$-adic valuation by $v_p$.
\ifkrenn\end{notation}\else\endrem\fi

\ifkrenn\begin{fact}\label{fact:alphabeta}\else
\def\alphabeta{Fact 2.1}
\pronounce{2.1 Fact}\fi
\item{\rm (1)} $\alpha_p(n)=v_p(n!)$. 
\item{\rm (2)} \ifdissertationkrenn For $1\le n\le p$, $\beta_p(n)=np$ and for
  $n>p$, $\beta_p(n)<np$. \else For $1\le k\le p$, $\beta_p(k)=kp$ and for
  $k>p$, $\beta_p(k)<kp$. \fi
\item{\rm (3)} For all $n\in\Z$, $v_p((n)_k)\ge \alpha_p(k)$; and
$v_p((k)_k)=v_p(k!)=\alpha_p(k)$.
\ifkrenn\end{fact}\else\endpronounce\fi

\ifkrenn\begin{proof}\else\proof\fi
Easy.
\ifkrenn\end{proof}\else\endproof\fi

\ifkrenn\begin{remark}\else\rem{Remark}\fi
The sequence $(\beta_p(n))_{n=1}^{\infty}$ is obtained by going
through the natural numbers in increasing order and repeating
each $k\in\N$ $v_p(k)$ times. 
For instance, 
$\beta_2(n)$ for $n\ge 1$ is: 
$2,4,4,6,8,8,8,10,12,12,14,16,16,16,16,18,20,20,\ldots$.
\ifkrenn\end{remark}\else\endrem\fi

The falling factorials $(x)_0=1$, $(x)_k=x(x-1)\ldots(x-k+1)$, $k> 0$,
form a basis of the free $\Z$-module $\Z[x]$, and representation with 
respect to this basis gives a convenient canonical form for a polynomial 
representing a given polynomial function on $\Zpn$.

\ifkrenn\begin{fact}[cf.~Keller and Olson {\rm \cite{KeOl68PF}}]\else
\pronounce{2.2 Fact}\(cf.~Keller and Olson {\rm \cite{KeOl68PF}}\)\ \fi
A polynomial $f\in\Z[x]$, $f=\sum_k a_k\; (x)_k$, induces the
zero-function mod $p^n$ if and only if 
$a_k\congr 0$ mod $p^{n-\alpha(k)}$
for all $k$ \(or, equivalently, for all $k<\beta(n)$\).
\ifkrenn\end{fact}\else\endpronounce\fi

\ifkrenn\begin{proof}\else\proof\fi
Induction on $k$ using the facts that $(m)_k=0$ for $m<k$, that
$v_p((n)_k)\ge \alpha_p(k)$ for all $n\in\Z$, and that
$v_p((k)_k)=v_p(k!)=\alpha_p(k)$.
\ifkrenn\end{proof}\else\endproof\fi

\ifkrenn\begin{corollary}[cf.~Keller and Olson {\rm \cite{KeOl68PF}}]\else
\pronounce{2.3 Corollary}\(cf.~Keller and Olson {\rm \cite{KeOl68PF}}\)\ \fi
Every polynomial function on $\Zpn$ is represented by a unique
$f\in\Z[x]$ of the form $f=\sum_{k=0}^{\beta(n)-1} a_k\; (x)_k$, 
with $0\le a_k<p^{n-\alpha(k)}$ for all $k$.
\ifkrenn\end{corollary}\else\endpronounce\fi

Comparing the canonical forms of polynomial functions mod $p^n$
with those mod $p^{n-1}$ we see that every polynomial function
mod $p^{n-1}$ gives rise to $p^{\beta(n)}$ different polynomial
functions mod $p^n$:

\ifkrenn\begin{corollary}[cf.~Keller and Olson {\rm \cite{KeOl68PF}}]
\label{cor:Fn}\else
\def\Fn{Corollary 2.4} \pronounce{2.4 Corollary}\(cf.~Keller and Olson {\rm
  \cite{KeOl68PF}}\)\ \fi Let $(F_n,\circ)$ be the monoid of polynomial
functions on $\Zpn$ with respect to composition and $\pi_{n}\colon
F_{n+1}\rightarrow F_{n}$ the canonical projection.
\item{\rm (1)}
For all $n\ge 1$ and for each $f\in F_{n}$ we have $|\pi_n^{-1}(f)|=p^{\beta(n+1)}$.
\item{\rm (2)}
For all $n\ge 1$, the number of polynomial functions on $\Zpn$ is 
$$|F_n|=p^{\smallsum_{k=1}^n \beta(k)}.$$
\ifkrenn\end{corollary}\else\endpronounce\fi

\ifdissertationkrenn Recall the following notation already given in the
introduction. \fi

\ifkrenn\begin{notation}\else\rem{Notation}\fi
We write $[f]_{p^n}$ for the function defined by $f\in\Z[x]$ on $\Zpn$.
\ifkrenn\end{notation}\else\endrem\fi


\ifkrenn\begin{lemma}\label{lem:carlitz}\else
\def\carlitz{Lemma 2.5}
\pronounce{2.5 Lemma}\fi
Every polynomial
$f\in\Z[x]$ is uniquely representable as
$$ f(x) = f_0(x) + f_1(x)(x^p-x) + f_2(x)(x^p-x)^2 + \ldots + f_m(x)(x^p-x)^m
+ \ldots$$
with $f_m\in\Z[x]$, $\deg f_m<p$, for all $m\ge 0$.
Now let  $f,g\in \Z[x]$.

\item{\rm (1)} If $n\le p$, then
$[f]_{p^n} = [g]_{p^n}$ is equivalent to:
$f_k = g_k$  mod $p^{n-k}\Z[x]$ for $0\le k< n$.
\item{\rm (2)}
$[f]_{p^2} = [g]_{p^2}$ is equivalent to:
$f_0 = g_0$ mod $p^2\Z[x]$ and $f_1=g_1$ mod $p\Z[x]$.
\item{\rm (3)}
$[f]_p = [g]_p$ and $[f']_p = [g']_p$
is equivalent to:
$f_0 = g_0$ mod $p\Z[x]$ and $f_1=g_1$ mod $p\Z[x]$.
\ifkrenn\end{lemma}\else\endpronounce\fi

\ifdissertationkrenn Note that (2) is just the special case of (1) with
$n=2$.\fi

\ifkrenn\begin{proof}\else\proof\fi
The canonical representation is obtained by repeated division with
remainder by $(x^p-x)$, and uniqueness follows from uniqueness of
quotient and remainder of polynomial division. Note that
$[f]_p = [f_0]_p$ and $[f']_p = [f_0'-f_1]_p$. This gives (3).

Denote by $f\sim g$ the equivalence relation 
$f_k = g_k$ mod $p^{n-k}\Z[x]$ for $0\le k< n$. Then
$f\sim g$ implies $[f]_{p^n} = [g]_{p^n}$. There are 
$p^{p+2p+3p+\ldots+np}$ equivalence classes of $\sim$
and $p^{\beta(1)+\beta(2)+\beta(3)+\ldots+\beta(n)}$
different $[f]_{p^n}$. For $k\le p$, $\beta(k)=kp$. 
Therefore the equivalence relations $f\sim g$ and
$[f]_{p^n} = [g]_{p^n}$ coincide. This gives (1)\ifdissertationkrenn.\else, and 
(2) is just the special case $n=2$.\fi
\ifkrenn\end{proof}\else\endproof\fi

We can rephrase this in terms of ideals of $\Zx$.

\ifkrenn\begin{corollary}\else
\pronounce{2.6 Corollary}\fi
For every $n\in\N$, consider the two ideals of $\Z[x]$
$$I_n=\{f\in\Z[x]\mid f(\Z)\subseteq p^n\Z\}
\qquad\hbox{\rm and}\qquad
J_n = (\{p^{n-k}(x^p-x)^k\mid 0\le k\le n\}).$$
Then
$[\Zx\colon I_n]=p^{\beta(1)+\beta(2)+\beta(3)+\ldots+\beta(n)}$
and
$[\Zx\colon J_n]=p^{p+2p+3p+\ldots+np}$.
Therefore, $J_n=I_n$ for $n\le p$, whereas for $n>p$, $J_n$ is properly
contained in $I_n$.
\ifkrenn\end{corollary}\else\endpronounce\fi

\ifkrenn\begin{proof}\else\proof\fi
$J_n\subseteq I_n$. The index of $J_n$ in $\Zx$ is
$p^{p+2p+3p+\ldots+np}$, because $f\in J_n$ if and only
if $f_k = 0$ mod $p^{n-k}\Z[x]$ for $0\le k< n$ in the 
canonical representation of \ifkrenn Lemma~\ref{lem:carlitz}\else\carlitz\fi.
The index of $I_n$ in $\Zx$
is $p^{\beta(1)+\beta(2)+\beta(3)\ldots+\beta(n)}$
by 
\ifkrenn Corollary~\ref{cor:Fn}\else\Fn\fi\ 
(2) and $[\Zx\colon I_n]<[\Zx\colon J_n]$
if and only if $n>p$ by 
\ifkrenn Fact~\ref{fact:alphabeta}\else\alphabeta\fi\ (2).
\ifkrenn\end{proof}\else\endproof\fi

\ifkrenn\begin{fact}[cf.~McDonald \cite{McD74FRI}]\label{fact:pp}\else
\def\ppfact{Fact 2.7}
\pronounce{2.7 Fact}\(cf.~McDonald {\rm \cite{McD74FRI}}\)\ \fi Let $n\ge 2$.
The function on $\Zpn$ induced by a polynomial $f\in\Z[x]$ is a
permutation if and only if
\item{\rm (1)} $f$ induces a permutation of $\Z_{p}$ and
\item{\rm (2)} the derivative $f'$ has no \ifdissertationkrenn root \else
  zero \fi mod $p$.  \ifkrenn\end{fact}\else\endpronounce\fi

\ifkrenn\begin{lemma}\label{lem:det}\else
\def\detlem{Lemma 2.8}
\pronounce{2.8 Lemma}\fi Let  $[f]_{p^n}$ and $[f]_p$ be the functions
defined by $f\in \Z[x]$ on $\Z_{p^n}$ and $\Zp$, respectively, and
$[f']_p$ the function defined by the formal derivative of $f$ on $\Zp$.
Then
\item{\rm (1)}
$[f]_{p^2}$ determines not just $[f]_p$, but also $[f']_p$.

\item{\rm (2)}
Let $n\ge 2$. Then $[f]_{p^n}$ is a permutation if and only if
$[f]_{p^2}$ is a permutation.

\item{\rm (3)}
For every pair of functions $(\alpha,\beta)$, 
$\alpha\colon\Zp\rightarrow \Zp$, $\beta\colon\Zp\rightarrow \Zp$,
there are exactly $p^p$ polynomial functions $[f]_{p^2}$ on $\Zptwo$ with
$[f]_p=\alpha$ and $[f']_p=\beta$.

\item{\rm (4)}
For every pair of functions $(\alpha,\beta)$, 
$\alpha\colon\Zp\rightarrow \Zp$ bijective,
$\beta\colon\Zp\rightarrow \Zp\setminus \{0\}$,
there are exactly $p^p$ polynomial permutations
$[f]_{p^2}$ on $\Zptwo$ with $[f]_p=\alpha$ and $[f']_p=\beta$.
\ifkrenn\end{lemma}\else\endpronounce\fi

\ifkrenn\begin{proof}\else\proof\fi
(1) and (3) follow immediately from 
\ifkrenn Lemma~\ref{lem:carlitz}\else\carlitz\fi\ for $n=2$ and 
(2) and (4) then follow from \ifkrenn Fact~\ref{fact:pp}\else\ppfact\fi.
\ifkrenn\end{proof}\else\endproof\fi

\ifkrenn\begin{remark}\else
\rem{2.9 Remark}\fi
\ifkrenn Fact~\ref{fact:pp} and Lemma~\ref{lem:det}%
\else\ppfact\ and \detlem\fi\ (2) imply that

\item{(1)} for all $n\ge 1$,
the image of $G_{n+1}$ under $\pi_n\colon F_{n+1}\rightarrow F_n$ is 
contained in $G_n$ and
\item{(2)} for all $n\ge 2$,
the inverse image of $G_n$ under
$\pi_n\colon F_{n+1}\rightarrow F_n$ is $G_{n+1}$.

We denote by $\pi_{n}\colon G_{n+1}\rightarrow G_n$ the restriction
of $\pi_{n}$ to $G_{n}$. This is the canonical epimorphism from the group
of polynomial permutations on $\Z_{p^{n+1}}$ onto the group of polynomial 
permutations on $\Zpn$.
\ifkrenn\end{remark}\else\endrem\fi

The above remark allows us to draw conclusions on the projective
system of groups $G_n$ from the information 
in \ifkrenn Corollary~\ref{cor:Fn}\else \Fn\ \fi concerning the projective
system of monoids $F_n$.

\ifkrenn\begin{corollary}\label{cor:inverseimage}\else
\def\inverseimage{Corollary 2.10}
\pronounce{2.10 Corollary}\fi
Let $n\ge 2$, and $\pi_{n}\colon G_{n+1}\rightarrow G_n$ the canonical
epimorphism from the group of polynomial permutations on 
$\Z_{p^{n+1}}$ onto the group of polynomial permutations on $\Zpn$.
Then $$|\ker(\pi_{n})|=p^{\beta(n+1)}.$$
\ifkrenn\end{corollary}\else\endpronounce\fi

\ifkrenn\begin{corollary}[cf.~Kempner {\rm \cite{Kem21PRS}} and 
Keller and Olson {\rm \cite{KeOl68PF}}]\else
\pronounce{2.11 Corollary}\(cf.~Kempner {\rm \cite{Kem21PRS}} and 
Keller and Olson {\rm \cite{KeOl68PF}}\)\ \fi
The number of polynomial permutations on $\Zptwo$ is
$$|G_2|=p!(p-1)^p p^p,$$
and for $n\ge 3$ the number of polynomial
permutations on $\Zptwo$ is 
$$|G_n|=p!(p-1)^p p^p p^{\smallsum_{k=3}^n \beta(k)}.$$
\ifkrenn\end{corollary}\else\endpronounce\fi

\ifkrenn\begin{proof}\else\proof\fi
In the canonical representation of $f\in\Zx$ in 
\ifkrenn Lemma~\ref{lem:carlitz}\else\carlitz\fi, there are
$p!(p-1)^p$ choices of coefficients mod $p$ for $f_0$ and $f_1$
such that the criteria of 
\ifkrenn Fact~\ref{fact:pp}\else\ppfact\fi\ for a polynomial permutation
on $\Zptwo$ are satisfied. And for each such choice there are
$p^p$ possibilities for the coefficients of $f_0$ mod $p^2$.
The coefficients of $f_0$ mod $p^2$ and those of $f_1$ mod $p$
then determine the polynomial function mod $p^2$.
So $|G_2|=p!(p-1)^p p^p$. The formula for $|G_n|$ then follows from
\ifkrenn Remark~\ref{rem:inverseimage}\else\inverseimage\fi.
\ifkrenn\end{proof}\else\endproof\fi

This concludes our review of polynomial functions and polynomial
permutations on $\Zpn$.
We will now introduce a homomorphic image of $G_2$ whose Sylow
$p$-groups bijectively correspond to the Sylow $p$-groups of $G_n$
for any $n\ge 2$.
\ifkrenn\else\goodbreak\fi

\ifkrenn\section{A group between $G_1$ and $G_2$}\else
\secthead{3.}{A group between $G_1$ and $G_2$}%
\nobreak\fi 
Into the projective system of monoids $(F_n,\circ)$ we insert an
extra monoid $E$ between $F_1$ and $F_2$ by means of monoid
epimorphisms $\theta\colon F_2\rightarrow E$ and 
$\psi\colon E\rightarrow F_1$ with $\psi\theta=\pi_1$. 

$$
F_1 \namedleftarrow\psi E\namedleftarrow\theta F_2\pileftarrow2
F_3\pileftarrow3\ldots
$$

The restrictions of $\theta$ to $G_2$ and of $\psi$ to the group of
units $H$ of $E$ will be group-epimorphisms, so that we also insert
an extra group $H$ between $G_1$ and $G_2$ into the projective
system of the $G_i$. 

$$
G_1 \namedleftarrow\psi H\namedleftarrow\theta G_2\pileftarrow2 
G_3\pileftarrow3\ldots
$$

In the following definition of $E$ and $H$, $f$ and $f'$ are just
two different names for functions. The connection with polynomials
and their formal derivatives suggested by the notation will appear
when we define $\theta$ and $\psi$.

\ifkrenn\begin{definition}\else\rem{Definition}\fi
We define the semi-group $(E,\circ)$ by
$$E=\{(f,f')\mid f\colon \Zp\rightarrow \Zp\;
f'\colon\Zp\rightarrow \Zp\}$$
(where $f$ and $f'$ are just symbols) with law of composition
$$(f,f')\circ (g,g') = (f\circ g,\; (f'\circ g)\cdot g'). $$
Here $(f\circ g)(x)= f(g(x))$ and 
$((f'\circ g)\cdot g')(x)=f'(g(x))\cdot g'(x)$.

We denote by $(H,\circ)$ the group of units of $E$.
\ifkrenn\end{definition}\else\endrem\fi

The following facts are easy to verify:

\ifkrenn\begin{lemma}\else
\pronounce{3.1 Lemma}\fi
\item{\rm (1)}
The identity element of $E$ is $(\id,1)$, with $\id$ denoting the
identity function on $\Zp$ and $1$ the constant function $1$.
\item{\rm (2)}
The group of units of $E$ has the form
$$H=\{(f,f')\mid f\colon \Zp\rightarrow \Zp\;\hbox{\rm bijective},\;
f'\colon\Zp\rightarrow \Zp\takeaway\{0\}\}.$$
\item{\rm (3)}
The inverse of $(g,g')\in H$ is 
$$(g,g')^{-1}=(g^{-1}, {1\over{g'\circ g^{-1}}}),$$
where $g^{-1}$ is the inverse permutation of the permutation $g$
and $1/a$ stands for the multiplicative inverse of a non-zero element
$a\in\Zp$, such that
$$({1\over{g'\circ g^{-1}}})(x)= {1\over{g'(g^{-1}(x))}}$$
means the multiplicative inverse in $\Zp\setminus\{0\}$ of $g'(g^{-1}(x))$.
\ifkrenn\end{lemma}\else\endpronounce\fi

Note that $H$ is a semidirect product of (as the normal subgroup)
a direct sum of $p$ copies of the cyclic group of order $p-1$
and (as the complement acting on it) the symmetric group on $p$ 
letters, $S_p$, acting on the direct sum by permuting its components. 
In combinatorics, one would call this a
wreath product (designed to act on the left) of the abstract group 
$C_{p-1}$ by the permutation group $S_p$ with its standard action on $p$
letters.  (Group theorists, however, have a narrower definition of wreath 
product, which is not applicable here.)

Now for the homomorphisms $\theta$ and $\psi$.

\ifkrenn\begin{definition}\else\rem{Definition}\fi
We define $\psi\colon E\longrightarrow F_1$ by $\psi(f,f')=f$. 
As for $\theta\colon F_2\rightarrow E$, given an element $[g]_{p^2}\in F_2$,
set $\theta([g]_{p^2})=([g]_p,[g']_p)$. $\theta$ is well-defined by
\ifkrenn Lemma~\ref{lem:det}\else\detlem\fi\ (1).
\ifkrenn\end{definition}\else\endrem\fi

\ifkrenn\begin{lemma}\label{lem:tlem}\else
\def\tlem{Lemma 3.2}
\pronounce{3.2 Lemma}\fi
\item{\(i\)}
$\theta\colon F_2\rightarrow E$ is a monoid-epimorphism.
\item{\(ii\)}
The inverse image of $H$
under $\theta\colon F_2\rightarrow E$ is $G_2$.
\item{\(iii\)}
The restriction of $\theta$ to $G_2$ is a group epimorphism
$\theta\colon G_2\rightarrow H$ with $|ker(\theta)|=p^p$.
\item{\(iv)}
$\psi\colon E\rightarrow F_1$ is a monoid epimorphism and
$\psi$ restricted to $H$ is a group-epimorphism
$\psi\colon H\rightarrow G_1$.
\ifkrenn\end{lemma}\else\endpronounce\fi

\ifkrenn\begin{proof}\else\proof\fi
(i) follows from \ifkrenn Lemma~\ref{lem:det}\else\detlem\fi\ (3)
and (ii) from \ifkrenn Fact~\ref{fact:pp}\else\ppfact\fi. 
(iii) follows from \ifkrenn Lemma~\ref{lem:det}\else\detlem\fi\ (4). 
Finally, (iv) holds
because every function on $\Zp$ is a polynomial function and
every permutation of $\Zp$ is  a polynomial permutation.
\ifkrenn\end{proof}\else\endproof\fi
\ifkrenn\else\goodbreak\fi

\ifkrenn\section{Sylow subgroups of $H$}\else
\secthead{4.}{Sylow subgroups of $H$}%
\nobreak\noindent\fi
We will first determine the Sylow $p$-groups of $H$. The Sylow $p$-groups
of $G_n$ for $n\ge 2$ are obtained in the next section as the inverse images 
of the Sylow $p$-groups of $H$ under the epimorphism $G_n\rightarrow H$.

\ifkrenn\begin{lemma}\label{lem:Slem}\else
\def\Slem{Lemma 4.1}
\pronounce{4.1 Lemma}\fi
Let $C_0$ be the subgroup of $S_p$ generated by the $p$-cycle 
$(0\,1\,2\ldots p-1)$. Then one Sylow $p$-subgroup of $H$ is
$$S= \{(f,f')\in H\mid f\in C_0,\; f'=1\},$$
where $f'=1$ means the constant function $1$.
The  normalizer of $S$ in $H$ is
$$N_H(S)= 
\{(g,g')\mid g\in N_{S_p}(C_0),\; g' \hbox{\rm\ a non-zero constant\ }\}.$$
\ifkrenn\end{lemma}\else\endpronounce\fi

\ifkrenn\begin{proof}\else\proof\fi
As $|H| = p! (p-1)^p$, and $S$ is a subgroup of $H$ of order $p$,
$S$ is a Sylow $p$-group of $H$.
Conjugation of $(f,f')\in S$ by $(g,g')\in H$ (using the fact that
$f'=1$) gives
$$(g,g')^{-1}(f,f')(g,g') = 
(g^{-1}, {1\over{g'\circ g^{-1}}})(f\circ g, g')=
(g^{-1}\circ f\circ g, 
{{g'}\over {g'\circ g^{-1}\circ f\circ g}})$$
The first coordinate of $(g,g')^{-1}(f,f')(g,g')$ being in $C_0$ for all
$(f,f')\in S$ is equivalent to $g\in N_{S_p}(C_0)$.
The second coordinate of $(g,g')^{-1}(f,f')(g,g')$ being the constant
function $1$ for all $(f,f')\in S$ is equivalent to
$$\forall x\in\Zp\quad  g'(x) = g'(g^{-1}(f(g(x))),$$
which is equivalent to $g'$ being constant on every cycle of $g^{-1}fg$,
which is equivalent to $g'$ being constant on $\Zp$, since $f$ can be 
chosen to be a $p$-cycle.
\ifkrenn\end{proof}\else\endproof\fi

\ifkrenn\begin{lemma}\else
\pronounce{4.2 Lemma}\fi
Another way of describing the normalizer of $S$ in $H$ is 
$$N_H(S)= \{(g,g')\in H\mid 
\exists k\ne 0\;\; \forall a,b\;\; g(a)-g(b)=k(a-b);\ 
\hbox{\rm $g'$ a non-zero constant}\}.$$
Therefore, $|N_H(S)|= p(p-1)^2$ and $[H\colon N_H(S)]=(p-1)!(p-1)^{p-2}$.
\ifkrenn\end{lemma}\else\endpronounce\fi

\ifkrenn\begin{proof}\else\proof\fi
Let $\sigma=(0\,1\,2\ldots p-1)$ and $g\in S_p$ then 
$$g\sigma g^{-1} = (g(0)\,g(1)\,g(2)\ldots g(p-1))$$
Now $g\in N_{S_p}(C_0)$ if and only if, for some $1\le k<p$
$g\sigma g^{-1} = \sigma^k$, i.e.,
$$(g(0)\;\; g(1)\;\; g(2)\; \ldots\; g(p-1))= (0\;\; k\;\; 2k\;\ldots\; (p-1)k),$$
all numbers taken mod $p$. This is equivalent to $g(x+1)=g(x)+k$ or
$$g(x+1)-g(x)=k$$ and further equivalent to $g(a)-g(b)=k(a-b)$.
Thus $k$ and $g(0)$ determine $g\in N_{S_p}(C_0)$, and there are
$(p-1)$ choices for $k$ and $p$ choices for $g(0)$. Together with the
$(p-1)$ choices for the non-zero constant $g'$ this makes $p(p-1)^2$
elements of $N_H(S)$.
\ifkrenn\end{proof}\else\endproof\fi

\ifkrenn\begin{corollary}\else
\def\HSylownum{Corollary 4.3}
\pronounce{4.3 Corollary}\fi 
There are $(p-1)!(p-1)^{p-2}$ Sylow $p$-subgroups of $H$.
\ifkrenn\end{corollary}\else\endpronounce\fi

\ifkrenn\begin{theorem}\label{thm:Hthm}\else
\def\Hthm{Theorem 4.4}
\pronounce{4.4 Theorem}\fi
The Sylow $p$-subgroups of $H$ are in bijective correspondence with
pairs $(C, \bar\varphi)$, where $C$ is a cyclic subgroup of order $p$
of $S_p$, $\varphi\colon\Zp\rightarrow \Zp\setminus\{0\}$ is a 
function and $\bar\varphi$ is the class of $\varphi$ with respect to the
equivalence relation of multiplication by a non-zero constant. The
subgroup corresponding to $(C, \bar\varphi)$ is
$$S_{(C, \bar\varphi)}= 
\{(f,f')\in H\mid f\in C,\; f'(x) = {{\varphi(f(x))}\over{\varphi(x)}}\}$$
\ifkrenn\end{theorem}\else\endpronounce\fi

\ifkrenn\begin{proof}\else\proof\fi
Observe that each $S_{(C, \bar\varphi)}$ is a subgroup of order $p$
of $H$. Different pairs $(C, \bar\varphi)$ give rise to
different groups: Suppose $S_{(C, \bar\varphi)}=S_{(D, \bar\psi)}$.
Then $C=D$ and for all $x\in\Zp$ and for all $f\in C$ we get
$${{\varphi(f(x))}\over{\varphi(x)}} = {{\psi(f(x))}\over{\psi(x)}}.$$
As $C$ is transitive on $\Zp$ the latter condition is equivalent to
$$\forall x,y\in \Zp\quad 
{{\psi(x)}\over{\varphi(x)}} = {{\psi(y)}\over{\varphi(y)}},$$
which means that $\varphi=k\psi$ for a nonzero $k\in\Zp$.

There are $(p-2)!$ cyclic subgroups of order $p$ of $S_p$,
and $(p-1)^{p-1}$ equivalence classes $\bar\varphi$ of functions
$\varphi\colon\Zp\rightarrow \Zp\setminus\{0\}$.  So the number of
pairs $(C, \bar\varphi)$ equals $(p-1)!(p-1)^{p-2}$, which is the
 number of Sylow $p$-groups of $H$, by the preceding corollary. 
\ifkrenn\end{proof}\else\endproof\fi

\ifkrenn\begin{proposition}\label{pro:intersect-sy}\else
\def\intpro{Proposition~4.5}
\pronounce{4.5 Proposition}\fi
If $p$ is an odd prime then the intersection of all Sylow $p$-subgroups
of $H$ is trivial, i.e., 
$$\bigcap_{(C,\overline\varphi)} S_{(C,\overline\varphi)}
= \{(\iota,1)\}.$$
If $p=2$ then $\left|H\right|=2$ and the intersection of all
Sylow $2$-subgroups of $H$ is $H$ itself.
\ifkrenn\end{proposition}\else\endpronounce\fi

\ifkrenn\begin{proof}\else\proof\fi
Let $p$ be an odd prime, and let $(f,f') \in \bigcap_{(C,\overline\varphi)}
S_{(C,\overline\varphi)}$. Suppose $f$ is not the identity function
and let $k\in\Zp$ such that $f(k)\neq k$.

Note that $\varphi$ in $(C,\overline\varphi)$ is arbitrary, apart from
the fact that $0$ is not in the image. Therefore, and because
$p\geq3$, among the various $\varphi$ there occur functions $\vartheta$ 
and $\eta$ with $\vartheta(k) = \eta(k)$ and 
$\vartheta(f(k)) \neq \eta(f(k))$. Now
$(f,f') \in S_{(D,\overline\vartheta)} \cap S_{(E,\overline\eta)}$ 
for any cyclic subgroups $D$ and $E$ of $S_p$ of order~$p$.

Therefore
$${{\vartheta(f(k))}\over{\vartheta(k)}} = f'(k) = 
{{\eta(f(k))}\over{\eta(k)}},$$
and hence $\vartheta(f(k)) = \eta(f(k))$, a contradiction.
Thus $f$ is the identity and therefore $f'=1$.

If $p=2$ then $\left|H\right|=2$ and therefore the one and only 
Sylow $2$-subgroup of $H$ is $H$.
\ifkrenn\end{proof}\else\endproof\fi

In the case $p\geq5$, the lemma above can be proved in a simpler way: There is
more than one cyclic group of order~$p$, so for $(f,f') \in
\bigcap_{(C,\overline\varphi)} S_{(C,\overline\varphi)}$, there are
distinct cyclic groups $D$ and $E$ of order~$p$ with $f\in D \cap E$.
Therefore $f$ has to be the identity.  

\ifkrenn\section{Sylow subgroups of $G_n$ and of the projective limit}\else
\secthead{5.}{Sylow subgroups of $G_n$ and of the projective limit $G$}%
\nobreak\noindent\fi

Again we consider the projective system of finite groups 
$$
G_1 \namedleftarrow\psi H\namedleftarrow\theta G_2\pileftarrow2
\ldots \pileftarrow{n-1} G_n \pileftarrow{n}
$$
where $(G_n,\circ)$ is the group of polynomial permutations on $\Z_{p^n}$ 
(with respect to composition of functions) and $H$ is the group
defined in section $3$. Let $G=\projlim G_n$ be the projective limit of
this system. Recall that a Sylow $p$-group of a pro-finite group is
defined as a maximal group consisting of elements whose order in
each of the finite groups in the projective system is a power of $p$.

\def\Sylownum{Theorem 5.1}
\pronounce{5.1 Theorem}
\item{\(i\)}
Let $(G_n,\circ)$ be the group of polynomial permutations 
on $\Z_{p^n}$ with respect to composition. 
If $n\ge 2$ there are $(p-1)!(p-1)^{p-2}$
Sylow $p$-groups of $G_n$. They are the inverse images of the
Sylow $p$-groups of $H$ (described in \Hthm) under the canonical
projection $\pi\colon G_n\rightarrow H$, with 
$\pi=\theta\pi_2\ldots\pi_{n-1}$.

\item{\(ii\)}
Let $G=\projlim G_n$. There are 
$(p-1)!(p-1)^{p-2}$ Sylow $p$-groups of $G$, which are the
inverse images of the Sylow $p$-groups of $H$ (described in \Hthm) 
under the canonical projection $\pi\colon G\rightarrow H$.
\endpronounce

\proof
In the projective system 
$G_1 \namedleftarrow\psi H\namedleftarrow\theta G_2\pileftarrow2
\ldots \pileftarrow{n-1} G_n$
the kernel of the group epimorphism $G_n\rightarrow H$ is a 
finite $p$-group for every $n\ge 2$, because
for $n\ge 2$ the kernel of $\pi_{n}\colon G_{n+1}\rightarrow G_{n}$
is of order $p^{\beta(n+1)}$ by 
\ifkrenn Remark~\ref{rem:inverseimage}\else\inverseimage\fi, 
and the kernel of
$\theta\colon G_2\rightarrow H$ is of order $p^p$ by
\ifkrenn Lemma~\ref{lem:tlem}\else\tlem\fi\ (iii).
So the Sylow $p$-groups of $G_n$ for $n\ge 2$ are
just the inverse images of the Sylow $p$-groups of $H$ and,
likewise, the Sylow $p$-groups of the projective limit $G$
are just the inverse images of the Sylow $p$-groups of $H$,
whose number was determined in \HSylownum.
\endproof

If we combine this information with the description of the
Sylow $p$-groups of $H$ in \Hthm\ we get the following explicit
description of the Sylow $p$-groups of $G_n$. Recall that $[f]_{p^n}$
denotes the function induced on $\Zpn$ by the polynomial
$f$ in $\Z[x]$ (or in $\Zpm[x]$ for some $m\ge n$).

\ifkrenn\begin{corollary}\else
\pronounce{5.2 Corollary}\fi
Let $n\ge 2$. Let $G_n$ be the group (with respect to composition) of
polynomial permutations on $\Z_{p^n}$. The Sylow $p$-groups of $G_n$ are
in bijective correspondence with
pairs $(C, \bar\varphi)$, where $C$ is a cyclic subgroup of order $p$
of $S_p$, $\varphi\colon\Zp\rightarrow \Zp\setminus\{0\}$ is a 
function and $\bar\varphi$ its class with respect to the
equivalence relation of multiplication by a non-zero constant. The
subgroup corresponding to $(C, \bar\varphi)$ is
$$S_{(C, \bar\varphi)}= 
\{[f]_{p^n}\in G_n\mid [f]_p\in C,\; 
[f']_p(x) = {{\varphi([f]_p(x))}\over{\varphi(x)}}\}.$$
\ifkrenn\end{theorem}\else\endpronounce\fi

\pronounce{Example}
A particularly easy to describe Sylow $p$-group of $G_n$ is the one
corresponding to $(C,\varphi)$ where $\varphi$ is a 
constant function and $C$ the subgroup of $S_p$ generated by
$(0\,1\,2\ldots p-1)$. It is the inverse image of $S$ defined in \ifkrenn
Lemma~\ref{lem:Slem}\else\Slem\fi\ and it consists of the 
functions on $\Zpn$ induced by polynomials $f$ such that the formal
derivative $f'$ induces the constant function $1$ on $\Zp$ and the 
function induced by $f$ itself on $\Zp$ is a power of $(0\,1\,2\ldots p-1)$.
\endpronounce

Combining \Sylownum\ with 
\ifkrenn Proposition~\ref{pro:intersect-sy}\else\intpro\fi\ 
we obtain the following description of the intersection of all 
Sylow $p$-groups of $G_n$ for odd $p$.

\ifkrenn\begin{corollary}\else
\pronounce{5.3 Corollary}\fi
Let $p$ be an odd prime.
\item{\(i\)} For $n\geq2$ the intersection of all Sylow $p$-groups of $G_n$ 
is the kernel of the projection $\pi\colon G\rightarrow H$.
\item{\(ii\)}
Likewise, the intersection of all Sylow $p$-groups of $G$ is the kernel
of the canonical epimorphism of $G$ onto $H$.
\item{\(iii\)}
The intersection of all Sylow $p$-groups of $G_n$ ($n\ge 2$) can also
be described as the normal subgroup
$$N = \{[f]_{p^n} \in G_n \mid [f]_p = \id, [f']_p=1\},$$ 
where $\id$ denotes the identity function on $\Zp$.
Its order is $p^p p^{\sum_{k=3}^n \beta(k)}$ and its index in $G_n$
(for $n\ge 2$) is $$[G_n:N] = p! (p-1)^p$$
\item{\(iv\)}
Likewise, the index of the intersection of all Sylow $p$-subgroups of
$G$ in $G$ is $p! (p-1)^p$.

\ifkrenn\end{corollary}\else\endpronounce\fi

\ifkrenn\begin{proof}\else\proof\fi
(i) and (ii) follow immediately from \Sylownum\ and \intpro. To see (iii),
let $\pi$ be the projection from $G_n$ to $H$ (that is $\pi =
\theta\pi_2\dots\pi_{n-1}$). Then $N$ is the inverse image of $\{(\iota,1)\}$,
the identity element of $H$, under $\pi$, and is therefore the intersection of
the Sylow $p$-groups of $G_n$ by (i). As the kernel of a group
homomorphism, $N$ is a normal subgroup.

The order of $N$ is the order of the kernel of $\pi$, which is the product
of $p^p$ (the order of the kernel of $\theta$) and $p^{\beta(k)}$ (the order
of the kernel of $\pi_{k-1}$) for $3\leq k \leq n$. Finally, the index of 
the kernel of the homomorphism of $G_n$ or $G$ onto $H$ is the order of $H$
which is $p!(p-1)^p$.
\ifkrenn\end{proof}\else\endproof\fi

\pronounce{Acknowledgement}
The authors wish to thank W.~Herfort for stimulating discussuions.
\endpronounce

\ifkrenn\else\line{References\hfill}\fi


\bibliography{ppsy}
\bibliographystyle{siam}

\bigskip

\line{
\vbox{
\hbox{S.\ F.\hfil}
\hbox{Institut f\"ur Mathematik A\hfil}
\hbox{Technische Universit\"at Graz\hfil}
\hbox{Steyrergasse 30\hfil}
\hbox{A-8010 Graz, Austria\hfil}
\tt
\hbox{frisch@tugraz.at\hfil}
}
\hskip1.5cm
\vbox{
\hbox{D.\ K.\hfil}
\hbox{Institut f\"ur Mathematik B\hfil\hfil}
\hbox{Technische Universit\"at Graz\hfil}
\hbox{Steyrergasse 30\hfil}
\hbox{A-8010 Graz, Austria\hfil}
\tt
\hbox{krenn@math.tugraz.at\hfil}
}
\hfill
}

\bye